\documentclass[10pt,a4]{article}

\usepackage{amsfonts}
\usepackage{amsmath}
\usepackage{amssymb}
\usepackage{amsthm}
\usepackage{amscd}
\usepackage{verbatim}
\usepackage[dvips]{graphicx}
\usepackage{graphics}
\usepackage[all,2cell]{xypic}

\theoremstyle{plain}
\newtheorem{thm}{Theorem}[section]
\newtheorem{prop}[thm]{Proposition}
\newtheorem{lem}[thm]{Lemma}
\newtheorem{cor}[thm]{Corollary}
\newtheorem{con}[thm]{Conjecture}
\newtheorem{parathm}[thm]{}
\newtheorem{mainthm}[thm]{Main Theorem}

\theoremstyle{definition}

\newtheorem{rem}[thm]{Remark}
\newtheorem{ex}[thm]{Example}

\theoremstyle{remark}
\newtheorem{para}[thm]{}

\DeclareMathOperator{\Cone}{Cone}
\DeclareMathOperator{\im}{Im}

\DeclareMathOperator{\Spec}{Spec}

\DeclareMathOperator{\coker}{Coker}
\DeclareMathOperator{\id}{id}
\DeclareMathOperator{\Ker}{Ker}
\DeclareMathOperator{\Hom}{Hom}
\DeclareMathOperator{\dmgm}{DM_{gm}}
\DeclareMathOperator{\dmgmeff}{DM_{gm}^{eff}}
\DeclareMathOperator{\dmeff}{DM^{eff}_{-}}
\DeclareMathOperator{\Homo}{H}

\DeclareMathOperator{\a1}{\mathbb{A}^1_k}
\DeclareMathOperator{\ztr}{\mathbb{Z}_{tr}}
\DeclareMathOperator{\ext}{Ext}
\DeclareMathOperator{\G_m}{\mathbb{G}_m}
\DeclareMathOperator{\ch_0}{CH_0}
\DeclareMathOperator{\Norm}{N}
\DeclareMathOperator{\mgm}{M_{gm}}
\DeclareMathOperator{\prodmgm}{Pro-DM_{gm}}
\DeclareMathOperator{\mgmc}{M_{gm}^c}
\DeclareMathOperator{\zequi}{\mathcal{Z}_{equi}}
\DeclareMathOperator{\smcor}{SmCor}
\DeclareMathOperator{\m}{M}

\DeclareMathOperator{\Pic}{Pic}
\DeclareMathOperator{\Shvnis}{Shv_{Nis}}
\DeclareMathOperator{\Cor}{Cor}
\DeclareMathOperator{\PST}{PST}
\DeclareMathOperator{\proa}{Pro-\mathcal{A}}
\DeclareMathOperator{\ord}{ord}
\DeclareMathOperator{\Gal}{Gal}
\DeclareMathOperator{\Jac}{Jac}

\DeclareMathOperator{\gl}{GL}
\DeclareMathOperator{\sm}{Sm}

\title{Motivic interpretation of Milnor $K$-groups attached to
Jacobian varieties}
\author{Satoshi Mochizuki
\footnote{This research is supported by the 21 century COE program at Graduate School of Mathematical Sciences, the University of Tokyo.
}}
\date{}

\begin{document}

\maketitle

\vspace{-1cm}
\begin{center}
\tt{mochi@ms.u-tokyo.ac.jp}
\end{center}

\renewcommand{\abstractname}{Abstract}
\begin{abstract}
In the paper \cite{Som90} p.105, Somekawa conjectures that his Milnor K-group $K(k,G_1,\ldots,G_r)$  attached to semi-abelian varieties $G_1$,\ldots,$G_r$ over a field $k$ is isomorphic to 
$\ext_{\mathcal{M}_k}^r(\mathbb{Z},G_1[-1] \otimes \ldots \otimes G_r[-1])$ where
$\mathcal{M}_k$ is a certain category of motives over $k$. The purpose of this note is to give remarks on this conjecture, when we take $\mathcal{M}_k$ as Voevodsky's category of motives $\dmeff(k)$ .
\end{abstract}


\tableofcontents

\section{Introduction}

To unify the Moore exact sequence and the Bloch exact sequence, K. Kato defined
the generalized Milnor $K$-groups attached to finite family of semi-abelian varieties over a base field $k$ in \cite{Som90}. (See also \cite{Akh00}, \cite{Kah92}.) That is, for semi-abelian varieties
$G_1,\ldots,G_r$ over $k$, he associated the group
$K(k,G_1,\ldots,G_r)$. (For precise definition, see \ref{Somekawa k group})
This group is a generalization of the Milnor $K$-group as the following example shows. 

\begin{ex}
In the notation above, if $G_1=G_2=\ldots=G_r=\G_m$, the following equality holds. 
$$K(k,\G_m,\ldots,\G_m)=K^M_r(k)$$
\end{ex}

\noindent
On the other hand this group is also a generalization of the Bloch group $V$. 

\begin{ex} \label{Bloch group}
Let $C$ be a projective smooth curve over $k$ such that $C(k) \ne \phi$.
We have the following equality
$$K(k,\Jac C,\G_m)=V(C).$$
where $V(C)$ is defined by S. Bloch (c.f. \cite{Blo81}) as the following way
$$V(C)=\frac{\scriptstyle {\Ker(\underset{x \in C^1}{\bigoplus} k(x)^\times 
\overset{\Sigma \Norm_{k(x)/k}}{\to} k^\times)}}
{\scriptstyle {\im(K_2(K(C)) 
\overset{\bigoplus \partial_x}{\to} \underset{x \in C^1}{\bigoplus} k(x)^\times)}}.$$
\end{ex}

\noindent
As is explained in \cite{Som90}, there is the generalized Bloch-Moore exact sequence.

\begin{para}
Let $k$ be a number field and $A$ a semi-abelian variety over $k$. We write
$G=\Gal(\bar{k}/k)$ and $G_v=\Gal(\bar{k_v}/k_v)$ for a place $v$. Let $T(A)$ be the Tate module of $A$. $S$ is a finite set of places including all 
Archimedian and places where $A$ has bad reduction.
Then $T(A)_G$ is a finite group by owing to \cite{KL81}. Let $m$ be a nonzero integer divisible by the order of $T(A)_G$. Somekawa proves the following generalized Moore-Bloch exact sequence (c.f. \cite{Som90} Theorem 4.1):
$$K(k,A,\G_m) \to \underset{v \notin S}{\bigoplus} T(A)_{G_v} \oplus
\underset{v \in S}{\bigoplus} K(k_v,A_v,\G_m)/m \to T(A)_G \to 0$$
In the case of $A=\G_m$, the above exact sequence is proved by Moore (c.f. \cite{Moo69})
$$K_2(k) \to \underset{v:\text{not complex}}{\bigoplus} \mu(k_v) \to \mu(k)
\to 0.$$
In the case of $A=\Jac C$ in the notation in \ref{Bloch group},
the above exact sequence is proved by S. Bloch, K. Kato and S. Saito
(c.f. \cite{Blo81}, \cite{KS83})
$$V(C) \to \underset{v \notin S}{\bigoplus} T(\Jac C)_{G_v} \oplus
\underset{v \in S}{\bigoplus} V(C \times_k k_v)/m \to T(\Jac C)_G \to 0.$$
\end{para}

\noindent
In \cite{Som90}, Somekawa conjectures that the Somekawa $K$-groups should be  
motivic cohomology groups attached to semi-abelian varieties. More precisely

\begin{con}(Somekawa conjecture) $ $
\\
Let  $G_1$,\ldots,$G_r$ be semi-abelian varieties over $k$, then
$K(k,G_1,\ldots,G_r)$ is isomorphic to  $\ext_{\mathcal{M}_k}^r(\mathbb{Z},G_1[-1] \otimes \ldots \otimes G_r[-1])$, where $\mathcal{M}_k$ is a certain category of motives over $k$ and $G_i[-1]$ means 1-motif (c.f \cite{Del74}).
\end{con}

\noindent
In this paper we will examine this conjecture, if we take $\mathcal{M}_k$
as Voevodsky's category of motives $\dmeff(k)$ .

\begin{mainthm} (Somekawa conjecture for Jacobian varieties) $ $
\\
Let $(C_1,a_1),\ldots,(C_n,a_n)$ be pointed projective smooth curves over perfect field $k$ which
admits resolution of singularities.  Then
$$K(k,\Jac C_1,\ldots,\Jac C_n) \overset{\sim}{\to} \Hom_{\dmeff(k)}
(\mgm (\Spec k), \mathbb{Z}(\underset{i=1}{\overset{n}{\wedge}} (C_i,a_1))[n]).$$
\end{mainthm}

\begin{para} \label{field assumption}
In this paper, let $k$ be a perfect field which admits resolution of singularity. 
\end{para}

\section{Milnor $K$-groups attached to semi-abelian varieties}

\subsection{Extension of valuations and tame symbols}

\begin{para}
Suppose $k$ is a field and $G$ is a semi-abelian variety defined over $k$,
that is, there is an exact sequence of group schemes (viewed as sheaves in
the flat topology) over $k$:
$$ 0 \to T \to G \to A \to 0$$
where $T$ is a torus and $A$ is an abelian variety.  
\end{para}

\begin{para}
In the notation above, let $K/k$ be an algebraic function field and $v$ 
a place of $K/k$. Let $L/K_v$ be a finite unramified Galois extension such that
$T \times_k F \overset{\sim}{\to} \G_m^n$ for the residue field $F$ of
$L$ and some $n$; let $w$ be the
unique extension of $v$ of $L$. We obtain the following commutative diagram
of exact sequences defining a map $r_w=(r_w^1,\ldots,r_w^n)$;\\
$$\xymatrix{
& 0 \ar[d] & 0 \ar[d]\\
0 \ar[r] & T(\mathcal{O}_w) \ar[r] \ar[d] & G(\mathcal{O}_w) \ar[r] \ar[d] 
& A(\mathcal{O}_w) \ar[r] \ar[d]^{\wr} & 0 \\
0 \ar[r] & T(L) \ar[r] \ar[d]_{\scriptstyle{\ord_w}} & G(L) \ar[r] \ar[d]^{\scriptstyle{r_w=(r_w^1,\ldots,r_w^n)}} & A(L)  \ar[r] & 0\\
  & \mathbb{Z}^n \ar[r]^{\scriptstyle{\ \ \ \id}} \ar[d] & \mathbb{Z}^n \ar[d]\\& 0 & 0          
}$$
\end{para}

\begin{para} \label{extended tame symbol}
In the notation above, we are going to construct a map
$$\partial_v:G(K_v) \otimes K_v^\times \to G(k(v)).$$
Fix $g \in G(K_v)$ and $h \in K_v^\times$. For each $i=1,\ldots,n$, we define
$h_i \in T(L)$ to be the $n$-th tuple having $h$ in the i-th coordinate and $1$
elsewhere. Then set
$$\varepsilon(g,h)=((-1)^{\ord_w(h)r^1_w(g)},\ldots,(-1)^{\ord_w(h)r^n_w(g)})
\in T(\mathcal{O}_w) \subset G(\mathcal{O}_w)$$
and
$$\tilde{\partial_v}(g,h)=\varepsilon(g,h)g^{\ord_w(h)}\prod_{i=1}^{n}
h_i^{-r^i_w(g)} \in G(\mathcal{O}_w).$$
We define the ``extended tame symbol'' $\partial_v(g,h)$ to be the image of
$\tilde{\partial_v}(g,h)$ under the canonical map 
$G(\mathcal{O}_w) \to G(F)$; Then $\partial_v(g,h)$
is invariant under the action of $\Gal(F/k(v))$, so that it belongs to 
$G(k(v))$. This definition of $\partial_v$ is independent of the choice of $L$ and of the isomorphic from the torus to
$\G_m^{\oplus n}$.
\end{para}

\subsection{Definition of the Milnor $K$-groups attached to semi-abelian varieties}

\begin{para} \label{Somekawa k group}
Let $k$ be a field and $G_1,\ldots,G_r$ a finite (possibly empty) family of 
semi-abelian varieties defined over $k$. We define Milnor $K$-groups
attached to semi-abelian varieties $K(k,G_1,\ldots,G_r)$ as follows.
If $r=0$, we write $K(k,\phi)$ for our groups and set $K(k,\phi)=\mathbb{Z}.$\\
For $r \geq 1$, we define
$$K(k,G_1,\ldots,G_r)=F/R$$
where
$$F=\underset{E/k:\text{finite}}{\bigoplus} G_1(E) \otimes \ldots \otimes G_r(E)$$and $R \subset F$ is the subgroup generated by the relation 
\textbf{R1}-\textbf{R2} below. \\
\textbf{R1} For any finite extensions $k \hookrightarrow  E_1
\overset{\phi}{\hookrightarrow} E_2$, let $g_{i_{0}} \in
G_{i_0}(E_2)$ and $g_i \in G_i(E_1)$ for $i \ne i_0$, the relation
$$(\phi^\ast(g_1) \otimes \ldots \otimes g_{i_0} \otimes \ldots \otimes 
\phi^\ast(g_r))_{E_2}- (g_1 \otimes \ldots \otimes \Norm_{E_2/E_1}(g_{i_0}) 
\otimes \ldots \otimes g_r)_{E_1}$$
(Here $\Norm_{E_2/E_1}$ denotes the norm map on the group scheme $G_{i_0}$)\\
\textbf{R2} For every algebraic function field $K/k$ and all choices
$g_i \in G_i(K)$, $h \in K^{\times}$ such that for each place $v$ of $K/k$, 
there exists $i(v)$ such that $g_i \in G_i(\mathcal{O}_v)$ 
for all $i \ne i(v)$, the relation
$$\sum_{v:\text{place of }K/k} (g_1(v) \otimes \ldots \otimes
\partial_v(g_{i(v)},h) \otimes \ldots \otimes g_r(v))_{k(v)/k}$$
Here $g_i(v) \in G_i(k(v))$ $(i \ne i(v))$ denotes the reduction of
$g_i \in \mathcal{O}_v$ modulo $m_v$ and $\partial_v(g_{i(v)},h)$ is the
extended tame symbol as defined in \ref{extended tame symbol}\\  
The class in $F/R$ of an element 
$a_1 \otimes \ldots \otimes a_r \in G_1(E) \otimes \ldots \otimes G_r(E)$
will be denoted
$\{a_1,\ldots,a_r\}_{E/k}$.
\end{para}

\begin{rem}
By the relation \textbf{R1}, if $\phi$ is a $k$-isomorphism $E_1 \overset{\sim}
{\to} E_2$, then 
$\{g_1,\ldots,g_r\}_{E_1/k}=\{\phi^\ast(g_1),\ldots,\phi^\ast(g_r)\}_{E_2/k}$.
This shows that symbols form a set.
\end{rem} 

\section{Triangulated categories of motives}

In this section, we will briefly review the definition of the triangulated
categories of motives. (c.f. \cite{TriCa}). 

\subsection{Triangulated category of effective geometric motives}

First we will review the construction of the category of geometric motives.

\begin{para}
1. Let $\sm/k$ be the category of schemes which are separated smooth, and of finite type over $k$. \\
2. Recall the definition of the category $\smcor(k)$: its objects are those of $\sm/k$. The set of morphism from $Y$ to $X$ is given by the group 
$\Cor(Y,X)$ of finite correspondences from $Y$ to $X$, defined as the free abelian group on the symbol $(Z)$, where $Z$ runs through the integral closed subschemes of $Y \times_k X$ which are finite over $Y$ and surjective over a connected component of $Y$.  We will denote the object of $\smcor(k)$ which corresponds to a smooth scheme $X$ by $[X]$.\\
\end{para}

\begin{para}
The category $\smcor(k)$ is an additive category. 
Consider the homotopy category
$\mathcal{H}^b(\smcor(k))$ of bounded complexes over $\smcor(k)$. 
Let $\mathcal{T}$ be the class of complexes of the following two forms:\\
1. For any smooth scheme $X$ over $k$ the complex
$$ [X \times \a1] \overset{[pr_1]}{\to} [X]$$
belongs to $\mathcal{T}$.\\
2. For any smooth scheme $X$ over $k$ and an open covering $X=U \cup V$ of $X$
the complex
$$[U \cap V] \overset{[j_U] \oplus [j_V]}{\to} [U] \oplus [V]
\overset{[i_U] \oplus (-[i_V])}{\to} [X]$$
belongs to $\mathcal{T}$. 
(here $j_U$, $j_V$, $i_U$, $i_V$ are the obvious open embeddings.)\\
Denote by $\bar{\mathcal{T}}$ the minimal thick subcategory of 
$\mathcal{H}^b(\smcor(k))$ which contains $\mathcal{T}$.\\
The triangulated category $\dmgmeff(k)$ of effective geometric motives over $k$
is the pseudo-Abelian envelope of the localization of 
$\mathcal{H}^b(\smcor(k))$ with respect to the thick subcategory 
$\bar{\mathcal{T}}$. We denote the obvious functor $Sm/k \to \dmgmeff(k)$ by
$\mgm$. 
\end{para}

\begin{para}
For a pair of smooth schemes $X$, $Y$ over $k$, we set
$$[X] \otimes [Y] := [X \times Y].$$
For any smooth schemes $X_1$, $Y_1$, $X_2$, $Y_2$ the external product of cycles defines a homomorphism:
$$\Cor(X_1,Y_1) \otimes \Cor(X_2,Y_2) \to \Cor(X_1 \times X_2, Y_1 \times Y_2)$$
which gives us a definition of tensor product of morphisms in $\smcor(k)$.
This structure defines in the usual way a tensor category structure on
$\mathcal{H}^b(\smcor(k))$ which can be descended to the category $\dmgmeff(k)$
by the universal property of localization.\\
Note that the unit object our tensor structure is $\mgm(\Spec k)$. We will denote it by $\mathbb{Z}$. 
\end{para}

\begin{ex} \label{indep rat pt}
Let $x,y:\Spec k \to \mathbb{P}^1_{k}$ be two $k$-rational points. Then
$\mgm(x)=\mgm(y):\mgm(\Spec k) \to \mgm(\mathbb{P}^1_{k})$.
\begin{proof}
We take an affine open set $\a1 \overset{j}{\hookrightarrow} \mathbb{P}^1_{k}$ which contains $x$ and $y$. That is, there are $\tilde{x},\tilde{y}:\Spec k \to
\a1$ such that $x=j \circ \tilde{x}$ and $y=j \circ \tilde{y}$. Then we have
$\mgm(x)=\mgm(j) \circ \mgm(\tilde{x})=\mgm(j) \circ \mgm(p)^{-1}
 =\mgm(j) \circ \mgm(\tilde{y})=\mgm(y)$, where $p:\a1 \to \Spec k$ is the structure morphism.
\end{proof} 
\end{ex}

\subsection{Triangulated category of effective motivic complexes}

To study the fundamental property of $\dmgmeff(k)$, 
Voevodsky uses sheaf theoretic method in \cite {TriCa}. 
More precisely, he constructs another category $\dmeff(k)$ using a sheaf category and he proves $\dmgmeff(k)$ admits a natural full embedding as a tensor category and a triangulated category to the category  $\dmeff(k)$.
We will review the construction of $\dmeff(k)$. 

\begin{para}
1. A presheaf with transfers on $\sm/k$ is an additive contravariant functor
from the category $\smcor(k)$ to the category of abelian groups.
We denote by $\PST(k)$ the category of presheaf with transfers on $\sm/k$.\\
2. A presheaf with transfers on $\sm/k$ is called a Nisnevich sheaf with 
transfers if the corresponding presheaf of abelian groups on 
$Sm/k$ is a sheaf in the Nisnevich topology.
We denoted by $\Shvnis (\smcor (k))$ the category of Nisnevich sheaves with transfers.
\end{para}

\begin{ex}
For any smooth scheme $X$ over $k$, a presheaf $\ztr(X):=\Cor(?,X)$ is 
a Nisnevich sheaf with transfers on $\sm/k$. (c.f. \cite{TriCa} Lemma 3.1.2).\\
For a $k$-rational point $x:\Spec k \to X$, we put 
$$\ztr(X,x):=\coker(
\ztr(\Spec k) \overset{\scriptstyle{\ztr(x)}}{\to} \ztr(X)).$$ 
\end{ex}

\begin{para}
1. A presheaf with transfers $F$ is called homotopy invariant if for any smooth
scheme $X$ over $k$ the projection $X \times \a1 \to X$ induces the isomorphism
$F(X) \to F(X \times \a1)$.\\
2. A Nisnevich sheaf with transfers is called homotopy invariant if it is
homotopy invariant as a presheaf with transfers.
\end{para}

\begin{para}
$\Shvnis (\smcor (k))$ is an abelian category. 
(c.f. \cite{TriCa} Theorem 3.1.4)
Inside the derived category $D^{-}(\Shvnis (\smcor (k)))$ of complexes bounded
from above, one defines the full subcategory $\dmeff(k)$ of effective motivic
complexes over $k$ as the one consisting of objects whose cohomology sheaves are homotopy invariant. This subcategory is triangulated. (Need the assumption of
perfectness of $k$.) (c.f. \cite{TriCa} Proposition 3.1.13).
\end{para}

\begin{para}
1. Let $F$ be a presheaf with transfers. There is a canonical surjection of presheaves
$$\bigoplus_{(X,x \in F(X))} \ztr(X) \to F.$$
Iterating this construction we get a canonical left resolution 
$\mathcal{L}(F)$ of $F$ which consists of 
direct sums of presheaves of the form $\ztr(X)$ for 
smooth schemes $X$ over $k$.\\
2. We set for two smooth schemes $X$,$Y$:
$$ \ztr(X) \otimes \ztr(Y) := \ztr(X \times Y)$$
and for two presheaves with transfers $F$, $G$:
$$F \otimes G := \underline{\Homo}_0(\mathcal{L}(F) \otimes \mathcal{L}(G)).$$ 
3. This construction provides us with a tensor structure on the derived 
category $D^{-}(\Shvnis(\smcor(k)))$. 
\end{para}

\noindent
To define the tensor structure on $\dmeff(k)$, we will need an alternative description of $\dmeff(k)$.

\begin{para}
Let $\Delta^\bullet$ be the standard cosimplicial object in $Sm/k$. For any presheaf with transfers $F$ on $Sm/k$ let $C_\ast(F)$ be the complex of presheaves on $Sm/k$ of the form $C_n(F)=\underline{\Hom}(\Delta^n,F)$ with differentials given by alternated sums of morphisms which correspond to the boundary morphisms
of $\Delta^\bullet$. This complex is called the singular simplicial complex of
$F$.\\
The following properties are fundamental.\\
1. If $F$ is a presheaf with transfers (resp. a Nisnevich sheaf with transfers)
then $C_\ast(F)$ is a complex of presheaves with transfers 
(resp. Nisnevich sheaves with transfers).\\
2. For any presheaf with transfers $F$ over $k$, the cohomology presheaf $h_i(F)$ of the complex $C_\ast(F)$ and its Nisnevich sheafication $h_i^{Nis}(F)$ are
homotopy invariant. (Need the assumption of perfectness of $k$.)
(c.f. \cite{TriCa} Lemma 3.2.1).\\
3. In view of 1. and 2., $C_\ast(?)$ is a functor from the category of 
Nisnevich sheaves with transfers on $Sm/k$ to $\dmeff(k)$.
\end{para}

\begin{prop} (c.f. \cite{TriCa} Proposition 3.2.3)
\\
The functor $C_\ast(?)$ can be extended to a functor
$$\mathcal{R}C:D^{-}(\Shvnis(\smcor(k)) \to \dmeff(k)$$
which is left adjoint to the natural embedding. The functor $\mathcal{R}C$
identifies $\dmeff(k)$ with localization of $D^{-}(\Shvnis(\smcor(k))$ with
respect to the localizing subcategory  $\mathcal{A}$ 
generated by complexes of the form
$$\ztr(X \times \a1) \overset{\ztr(pr_1)}{\to} \ztr(X)$$
for smooth schemes $X$ over $k$.
\end{prop}

\begin{para}
1. In the notation above, $\mathcal{A}$ is a $\otimes$-ideal, that is,
for any object $T$ of $D^{-}(\Shvnis(\smcor(k))$ and an object $S$ of
$\mathcal{A}$ the object $T \otimes S$ belongs to $\mathcal{A}$.
(c.f. \cite{TriCa} Lemma 3.2.4). \\
2. We define tensor structure on $\dmeff(k)$ as the descent of the tensor structure on $D^{-}(\Shvnis(\smcor(k))$ with respect to the projector $\mathcal{R}C$.
Note that such a descent exists by the universal property of localization 
and 1.
\end{para}

\begin{thm}(c.f. \cite{TriCa} Theorem 3.2.6) \label{embedding thm}
\\
There is a commutative diagram of functors of the form
such that the following conditions hold:\\
$$\xymatrix{
\mathcal{H}^b(\smcor(k)) \ar[r]^{\! \! \! L} \ar[d] & D^{-}(\Shvnis(\smcor(k)) \ar[d]^{\mathcal{R}C}\\
\dmgmeff(k) \ar[r]^{i} & \dmeff (k)
}$$
1. The functor $i$ is a full embedding with a dense image.\\
2. For any smooth scheme $X$ over $k$ the object $\mathcal{R}C(L(X))$ is
canonically isomorphic to the $C_\ast(\ztr(X))$.\\
3. All functors preserve tensor and triangulated structures.
\end{thm}

\begin{ex} \label{pointed motif}
Let $x:\Spec k \to X$ be a $k$-rational point of a smooth scheme $X$. Then
we have an identification
$$C_\ast(\ztr(X,x)) \overset{\sim}{\to} \Cone(\mgm(\Spec k) \overset{\scriptstyle{\mgm(x)}}{\to}
\mgm(X))$$ defined by
$$\begin{bmatrix} \xymatrix{ 
\ztr(\Spec k) \ar[d]_{\scriptstyle{\ztr(x)}}\\ \ztr(X)} \end{bmatrix}
\begin{matrix} \xymatrix{\to \ar@{}[d]\\ \to\\ } \end{matrix} 
\begin{bmatrix} \xymatrix{
0 \ar[d] \\ \ztr(X,x)} \end{bmatrix}.$$
\end{ex}

\subsection{Motives with compact support}

In this subsection, we will briefly review the notation and fundamental result of \cite{FV00}, \cite{RelCy} and \cite{TriCa}.  

\begin{para}
For any scheme of finite type $X$ over $k$ and any $r \geq 0$ we denote by
$\zequi (X,r)$ the presheaf on the category of smooth schemes over $k$ which
takes a smooth scheme Y to free abelian groups generated by closed integral subschemes $Z$ of $X \times Y$ which are equidimensional of relative dimension $r$ over $Y$.\\
This presheaf has the following property.\\
1. $\zequi (X,r)$ is a sheaf in the Nisnevich topology.\\ 
2. It has a canonical structure of a presheaf with transfers.\\  
3. The presheaf with transfers $\zequi (X,r)$ is covariantly functorial with respect to proper morphisms of $X$ by means of the usual proper push-forward of 
cycles. \\
4. It is contravariantly functorial with an appropriate dimension shift with respect to flat equidimensional morphisms.\\
5. There is a pairing 
$$ \times:\zequi(X,r) \otimes \zequi(X',r') \to \zequi(X \times X',r+r')$$
of presheaves.:\\
Let $U$ be a smooth scheme over $k$. For any pair of integral closed subschemes
$Z \subset X \times U, Z' \subset X' \times U$ equidimensional over $U$.
Sending $Z$, $Z'$ to the cycle associated to the subscheme 
$Z \times_U Z' \subset X \times X' \times U$ determines a pairing.\\
6. (Need the assumption of \ref{field assumption}.) The flat pull-back morphism induces a quasi-isomorphism
$$C_\ast(\zequi(X,0)) \to C_\ast(\zequi(X \times \mathbb{A}^n,n)).$$
\end{para}

\begin{para}
For any scheme of finite type over $k$ the object $C_\ast(\zequi (X,0))$ belongs to $\dmeff (k)$. Moreover it belongs to $\dmgmeff (k)$. (Using the assumption
\ref{field assumption}) (c.f. \cite{TriCa}, Corollary 4.1.6). 
We will denote $\mgmc (X):=C_\ast(\zequi (X,0))$ and
call it a motivic complex of $X$ with compact support.\\
Since $\ztr(X)$ is a subsheaf of $\zequi(X,0)$, the inclusion morphism induces the natural morphism $\mgm(X) \to \mgmc(X)$.\\ 
The following properties are fundamental. 
(c.f. \cite{TriCa}, Proposition 4.1.5, Proposition 4.1.7)\\
1. If $X$ is proper then the canonical morphism $\mgm (X) \to \mgmc (X)$
is the isomorphism.\\
2. (Need the assumption \ref{field assumption}.)  Let $Z$ be a closed subscheme
of $X$. Then there is a canonical distinguished triangle of the form
$$ \mgmc (Z) \to \mgmc (X) \to \mgmc (X-Z) \to \mgmc(Z)[1].$$  
3. For $X$, $Y$ of finite type over $k$, the pairing
$\zequi(X,0) \otimes \zequi(Y,0) \to \zequi(X \times_k Y,0)$
induces an isomorphism
$$\mgmc(X) \otimes \mgmc(Y) \overset{\sim}{\to} \mgmc(X \times_k Y).$$
\end{para}

\subsection{Tate object}

\begin{para}
For any smooth scheme $X$ over $k$, the morphism $X \to \Spec k$ gives us a morphism in $\dmgmeff(k)$ of the form $\mgm(X) \to\mathbb{Z}$. There is a canonical
split distinguished triangle
$$\widetilde{\mgm}(X) \to \mgm(X) \to \mathbb{Z} \to \widetilde{\mgm}(X)[1]$$
where $\widetilde{\mgm}(X)$ is the reduced motif of $X$ represented in
$\mathcal{H}^b(\smcor(k))$ by the complex $[X] \to [\Spec k]$.
\end{para}

\begin{ex} \label{cano isom}
In the notation above, for any $k$-rational point $x:\Spec k \to X$, we have
the canonical identification $\mgm(X,x) \overset{\sim}{\to} \widetilde{\mgm}(X)$
as the following way.
$$\begin{bmatrix} \xymatrix{
[\Spec k] \ar[d]^{\scriptstyle{x}} \\ [X] \ar[d] \\ 0 }\ \ \end{bmatrix}
\begin{matrix} 
\xymatrix{
\longrightarrow \ar@{}[d] \\ \overset{\scriptstyle{\id - \mgm(x \circ p)}}{\longrightarrow} \ar@{}[d] \\ \longrightarrow  }
\end{matrix}
\begin{bmatrix}
\xymatrix{
0 \ar[d] \\ [X] \ar[d]^{\scriptstyle{p}} \\ [\Spec k] } \end{bmatrix}$$
where $p:X \to \Spec k$ is the structure morphism.\\
$x:\Spec k \to X$ defines splitting $\mgm(X) \overset{\sim}{\to} \mgm(X,x) 
\oplus \mathbb{Z}$.
\end{ex}

\begin{para}
We define the Tate object $\mathbb{Z}(1)$ of $\dmgmeff(k)$ as 
$\widetilde{\mgm}(\mathbb{P}^1)[-2]$.
We further define $\mathbb{Z}(n)$ to be
the $n$-th tensor power of $\mathbb{Z}(1)$.\\
For any object $A$ of $\dmgmeff(k)$ we put  
$$A(n)= A \otimes \mathbb{Z}(n)$$
$$A\{n\}= A \otimes \mathbb{Z}(n)[n]$$
$$A((n))= A \otimes \mathbb{Z}(n)[2n].$$
By Example \ref{cano isom}, for any
$x:\Spec k \to \mathbb{P}^1_{k}$, we have the canonical isomorphism
$\mgm(\mathbb{P}^1_{k},x) \overset{\sim}{\to} \mathbb{Z}((1))$.
\end{para}

\begin{para} \label{Tate object}
Let $x:\Spec k \to \mathbb{P}^1_{k}$ be a $k$-rational point.
Comparing the following split distinguished triangles 
$$\xymatrix{
\mgmc(\Spec k) \ar[r] \ar[d]_{\scriptstyle{\id}} & \mgmc(\mathbb{P}^1) \ar[r] \ar[d]_{\scriptstyle{\id}}& \mgmc(\mathbb{A}^1) 
\ar[r] \ar[d] & \mgmc(\Spec k)[1] \ar[d]_{\scriptstyle{\id}}\\
\mgm(\Spec k) \ar[r]_{\scriptstyle{\mgm(x)}} & \mgm(\mathbb{P}^1_{k}) \ar[r] & \mathbb{Z}((1)) \ar[r] & \mgm(\Spec k)[1], 
}$$
where $\mgm(\mathbb{P}^1_{k}) \to \mathbb{Z}((1))$ is defined by 
$$\mgm(\mathbb{P}^1_{k}) \underset{\text{can}}{\to} \mgm(\mathbb{P}^1,x) 
\overset{\sim}{\to} \mathbb{Z}((1)),$$
we know that there is a natural isomorphism $\mathbb{Z}((1)) \overset{\sim}{\to}\mgmc(\mathbb{A}^1)$. It does not depend on the choice of a $k$-rational point by Example \ref{indep rat pt}.
Similarly using Mayer-Vietoris sequence for canonical covering of $\mathbb{P}^1_k$, we know also that there is a natural isomorphism
$\mathbb{Z}\{1\} \overset{\sim}{\to}
\widetilde{\mgm}(\mathbb{A}^1-\{0\})$.
\end{para}

\subsection{The triangulated category of geometric motives}

In this subsection, we will define the triangulated category $\dmgm(k)$ of geometric motives over $k$.

\begin{para}
1. We define the category $\dmgm(k)$: its objects are pairs of the form $(A,n)$ where $A$ is an object of $\dmgmeff(k)$ and $n \in \mathbb{Z}$ and morphisms are 
defined by the following formula 
$$\Hom_{\dmgm(k)}((A,n),(B,m)):=\lim_{k \geq -n,-m} 
\Hom_{\dmgmeff(k)}(A(k+n),B(k+m)).$$
2. The category $\dmgm(k)$ with the obvious shift functor and class of distinguished triangles is a triangulated category.\\
3. The permutation involution on $\mathbb{Z}(1) \otimes \mathbb{Z}(1)$
is identity in $\dmgmeff(k)$. (c.f. \cite{TriCa} Corollary 2.1.5)\\
4. Using the fact of 3. and general theory, $\dmgm(k)$ has a natural tensor 
structure.\\ 
\end{para}

\begin{thm}(c.f. \cite{Voe02} The cancellation theorem) \label{can thm}
\\
(Need the assumption of perfectness of $k$.) For objects
$A$, $B$ in $\dmgmeff(k)$ the natural map
$$? \otimes \id_{\mathbb{Z}(1)}:\Hom_{\dmgmeff(k)}(A,B) \to
\Hom_{\dmgmeff(k)}(A(1),B(1))$$
is an isomorphism. Thus the canonical functor
$$\dmgmeff(k) \to \dmgm(k)$$
is a full embedding.
\end{thm}

\section{Various morphisms between motives}

In this section, we will briefly review \cite{TriCa} and \cite{MotGe}

\subsection{Transpose for finite equidimensional morphisms}

\begin{para}
Let $X,Y$ be smooth schemes over $k$ and $f:X \to Y$ a finite equidimensional morphism. Then we have the transpose of $f$, ${}^tf:Y \to X$ in $\smcor(k)$. That is
$s(\Gamma_f) \in \Hom_{\smcor(k)}(Y,X)$, where $s:X \times Y \to Y \times X$ is the switch morphism.
\end{para}

\begin{ex} \label{norm}
Let $L/K/k$ be finite field extensions. Then we have a canonical morphism
$i:\Spec L \to \Spec K$. Hence we get 
$\mgm({}^ti):\mgm(\Spec K) \to \mgm(\Spec L)$ in $\dmgm(k)$. We shall write this map as
$\Norm_{L/K}$ because there is the following commutative diagram
$$\xymatrix{ 
\Hom (\mgm(\Spec L),\mathbb{Z}\{n\}) \ar[r]^{\sim} \ar[d]_{\scriptstyle{\Hom (\Norm_{L/K},\mathbb{Z}\{n\})}} & K^M_n(L) \ar[d]^{\scriptstyle{\Norm_{L/K}}}\\
\Hom (\mgm(\Spec K),\mathbb{Z}\{n\}) \ar[r]^{\sim} & K^M_n(K) .
}$$
(c.f. \cite{BKcon} Lemma 3.4.4. See also Theorem \ref{Milkisom})
\end{ex}

\subsection{Pull back for flat equidimensional morphisms}

\begin{para} \label{flat pull back def}
Let $X$, $Y$ be smooth schemes and
$f:X \to Y$ a flat equidimensional morphism of relative dimension $n$.
Then one can define a morphism 
$f^\ast:\mgmc(Y)((n)) \to \mgmc(X)$ as follows.
(This is slightly different from the definition in \cite{TriCa} Corollary 4.2.4).
\begin{eqnarray*}
\mgmc(Y)((n)) & = &  
C_\ast(\zequi(Y \times \mathbb{A}^n,0))\\
& \overset{\scriptstyle{{C_\ast(f \times \id_{\mathbb{A}^n})^\ast}}}{\to} & C_\ast(\zequi(X \times \mathbb{A}^n,n))\\ 
& \underset{\text{qis}}{\leftarrow} & C_\ast(\zequi(X,0))=\mgmc(X)
\end{eqnarray*} 
where every quasi-isomorphisms are induced from flat pull backs.\\
\end{para}

\begin{para} \label{fund of flat pullback}
The following properties are easily proved.\\
1 In the notation above, if $f$ is an open immersion, this morphism coincides with the canonical morphism $\mgmc (Y) \to \mgmc (X)$.\\
2 In the notation above, if  $X$ and $Y$ are proper over $\Spec k$ and $f$ is flat finite equidimensional, then $\mgm({}^t f)=f^\ast$.
(c.f. \cite{MotGe} Lemma 1.1.2)\\ 
3 Let $X$, $Y$ and $Z$ be smooth schemes and
$f:X \to Y$, $g:Y \to Z$ flat equidimensional morphisms of relative dimension $n$ and $m$ respectively. Then we have $f^\ast \circ g^\ast((m))=(g \circ f)^\ast$.
\end{para}

\begin{lem} \label{flat pillback example} $ $
\\
Let $p:\a1 \to \Spec k$ be the structure morphism. Then 
$p^\ast:\mgmc(\Spec k)((1))=\mgmc(\a1) \to \mgmc(\a1)$ coincides with 
$\id_{\mgmc(\a1)}$.
\end{lem}

\begin{proof} 
By definition (See \ref{flat pull back def}) , 
$$p^\ast=C_\ast(\zequi(\a1,0)) \overset{\scriptstyle{C_\ast(pr_1^\ast)}}{\to}
C_\ast(\zequi(\mathbb{A}^2_{k},1))
 \overset{\scriptstyle{(C_\ast(pr_2^\ast))^{-1}}}{\to} C_\ast(\zequi(\a1,0))$$ 
where $pr_1,pr_2$ are two projections $pr_1, pr_2:\mathbb{A}^2_k \to \a1$.
We assert that two projections $pr_1,pr_2$ induce the same morphism $pr_1^\ast=pr_2^\ast:C_\ast(\zequi(\a1,0)) \to C_\ast(\zequi(\mathbb{A}^2_k,1))$ in $\dmeff (k)$. Since
$$pr_1 =\begin{pmatrix} 0 &1\\ 1 &0 \end{pmatrix} \circ pr_2,$$
it suffices to prove that the action of $\gl_2(k)$ on $\mathbb{A}^2_{k}$ induces trivial action on $C_\ast(\zequi(\mathbb{A}^2,1))$ in $\dmeff(k)$ by flat pull back.
On the other hand $\gl_2(k)$ is generated by the elements of conjugate of
$\begin{pmatrix} 1 & 0\\ \ast & \ast \\ \end{pmatrix}$.
For 
$A=\begin{pmatrix} 1 & 0\\ c & d\\ \end{pmatrix} \in \gl_2(k)$, considering the
following diagram
$$\xymatrix{
C_\ast(\zequi(\mathbb{A}^2_{k},1)) \ar[r]^{A^\ast} & C_\ast(\zequi(\mathbb{A}^2_{k},1))\\
C_\ast(\zequi(\a1,0)) \ar[r]_{\id} \ar[u]^{\scriptstyle{pr_1^\ast}}_\wr & C_\ast(\zequi(\a1,0)) \ar[u]_{\scriptstyle{pr_1^\ast}}^{\wr},\\
}$$ 
we get the result.   
\end{proof}

\subsection{Motives with closed support}

\begin{para}
We call $(X,Z)$ a closed pair if $X$ is a smooth scheme over $k$ and $Z$ is a closed subscheme. If $Z$ is smooth over $k$, we call $(X,Z)$ a smooth pair. 
We call a pair of morphisms of scheme $(f,g):(Y,T) \to (X,Z)$ a morphism of closed pair if a commutative square 
$\xymatrix{
T \ar[r] \ar[d]_g & X \ar[d]^f\\
Z \ar[r] & Y
}$ 
is a Cartesian square as underlying topological spaces.  
Such a morphisms called Cartesian (resp. excisive) if the diagram above is a Cartesian square (resp. $f$ is \'etale and $g_{red}$ is an isomorphism.)
\end{para}

\begin{para}
Let $X$ be a smooth scheme and $U$ its open subset. Then we define
$$\mgm(X/U):=C^\ast(\coker(\ztr(U) \to \ztr(X))).$$ 
For any closed pair $(X,Z)$,
we define relative motif associated $(X,Z)$ by 
$\m_Z(X):=\mgm(X/X-Z)$. 
By definition there is a canonical distinguished triangle of the form
$$\mgm(X-Z) 
\overset{\mgm(j)}{\to} \mgm(X) 
\overset{i^\sharp}{\to} \m_Z(X) \to \mgm(X-Z)[1].$$
where $i:Z \hookrightarrow X$ is a closed immersion and 
$j:X-Z \hookrightarrow X$ is an open immersion. 
\end{para}

\begin{para}
For any morphisms of closed pair $(f,g):(Y,T) \to (X,Z)$, we associate a
morphism $(f,g)_\ast:\m_T (Y) \to \m_Z (X)$ which makes the following diagram
 commute
$$\xymatrix{
0 \ar[r] & \mgm(Y-T) \ar[r] \ar[d]_{\scriptstyle{\mgm (h)}} & \mgm(Y) \ar[r] \ar[d]_{\scriptstyle{\mgm (f)}}  & \m_T(Y) \ar[r] \ar[d]^{\scriptstyle{(f,g)_\ast}} & 0\\
0 \ar[r] & \mgm(X-Z) \ar[r] &  \mgm(X) \ar[r] & \m_Z(X) \ar[r] & 0
}$$
where $h:Y-T \to X-Z$ is a induced morphism from $f$.
\end{para}

\begin{para}
In the notation above, if $f$ is finite equidimensional, we associate a
morphism $(f,g)^\ast:\m_T (Y) \to \m_Z (X)$ which makes the following diagram
 commute
$$\xymatrix{
0 \ar[r]  & \mgm(Y-T) \ar[r] & \mgm(Y) \ar[r] & \m_T(Y) \ar[r] & 0\\
0 \ar[r] & \mgm(X-Z) \ar[r] \ar[u]^{\scriptstyle{\mgm ({}^t h)}} & \mgm(X) \ar[r] \ar[u]^{\scriptstyle{\mgm ({}^t f)}} & \m_Z(X) \ar[r] \ar[u]_{\scriptstyle{(f,g)^\ast}} & 0.
}$$
\end{para}

\begin{prop}(c.f. \cite{IntMo} Proposition 2.3)
\\
\textbf{(Red)} Reduction: If $(X,Z)$ is a closed pair, the canonical morphism
$(X,Z_{red}) \to (X,Z)$ induces identity map $\m_{Z_{red}}(X) \to \m_Z (X)$.\\
\textbf{(Add)} Additivity: Let $X$ be a smooth scheme, and $Z$, $W$ disjoint closed subschemes of $X$. Then induced morphism $\m_{Z \coprod W}(X) \to \m_Z (X) \oplus \m_W(X)$ is an isomorphism. \\
\textbf{(Exc)} Excision: Any excisive morphism $(Y,T) \to (X,Z)$ induces an isomorphism $\m_T (Y) \to \m_Z (X)$.\\
\textbf{(MV)} Mayer-Vietoris: Let $X$ be a smooth scheme over $k$, $U$ and $V$ two open subsets of $X$ such that $X=U \cup V$, and $Z$ a closed subscheme of 
$X$. Then we have a distinguished triangle of the form
$$\m_{Z \cap U \cap V}(U \cap V) \to \m_{Z \cap U}(U) \oplus \m_{Z \cap V}(V)
\to \m_Z(X) \overset{+1}{\to}.$$
\textbf{(Htp)} Homotopy invariance: A Cartesian morphism 
$\pi:(\mathbb{A}_X^1,\mathbb{A}_Z^1) \to (X,Z)$ induced from the canonical
projection induces an isomorphism 
$$\mgm(\pi):\m_{\mathbb{A}_Z^1}(\mathbb{A}_X^1) \to \m_Z (X).$$
\end{prop}

\subsection{Thom isomorphism}

\begin{para}
Let $X$ be a scheme and $E/X$ a vector bundle.
We consider $X$ as a closed subscheme of $E$ by zero section. We define the
Thom motif of $E/X$ by $\mgm (ThE):=\m_X(E)$. 
\end{para}

\begin{para}
In the notation above, if rank of $E$ is $n$, 
there is the Thom isomorphism 
$\theta(E):\mgm (ThE) \overset{\sim}{\to} \mgm(X)((n))$.
We will briefly review the construction of this isomorphism.
\end{para}

\begin{para}
Let $X$ be a smooth scheme and $\Delta_X:X \to X \times_k X$
the diagonal immersion. Let $\mathcal{M}$, $\mathcal{N}$ be objects of
$\dmeff(k)$ and
$$\alpha:\mgm(X) \to \mathcal{M}$$
$$\beta:\mgm(X) \to \mathcal{N}$$
are morphisms in $\dmeff(k)$.\\
We define external cup product of $\alpha$ and $\beta$ over $X$ is 
composition of the following morphisms 
$$\mgm(X) \overset{\mgm(\Delta_X)}{\to} \mgm(X) \otimes \mgm(X)
\overset{\alpha \otimes \beta}{\to} \mathcal{M} \otimes \mathcal{N}.$$
We denote this morphism by $\alpha \boxtimes_X \beta$ or 
$\alpha \boxtimes \beta$
\end{para}

\begin{para}
In the notation above, if $\mathcal{M}=\mathbb{Z}((m))$ and
$\mathcal{N}=\mathbb{Z}((n))$, then we have the canonical isomorphism
$\varepsilon:\mathbb{Z}((m)) \otimes \mathbb{Z}((n)) 
\overset{\sim}{\to}\mathbb{Z}((m+n))$.\\
We put $\alpha \cup \beta:=\varepsilon \circ \alpha \boxtimes \beta$
and call it internal cup product of $\alpha$ and $\beta$.
\end{para}

\begin{para}
There is a natural isomorphism as a presheaf with transfers (c.f. \cite{MotGe}
Corollary 2.2.7)
$$c_1:\Pic(?) \to \Hom_{\dmgm(k)}(\mgm(?),\mathbb{Z}((1))).$$
For any smooth scheme $X$ and $\mathcal{L} \in \Pic(X)$, we will call 
$c_1(\mathcal{L})$ a motivic Chern class of $\mathcal{L}$. 
\end{para}

\begin{para}
Let $X$ be a smooth scheme, $E$ a vector bundle over $X$ of rank $n$,
$\lambda_E$ the canonical invertible sheaf on $\mathbb{P}(E)$ and
$p:\mathbb{P}(E) \to X$ the canonical projection. We define the $r$-th 
motivic Lefschetz projector of $E$ by 
$$l_r(E):=c_1(\lambda_E)^{\cup r} \boxtimes \mgm(p):\mgm(\mathbb{P}(E)) \to 
\mgm(X)((r))$$
We put the motivic Lefschetz operator as
$$l(E):= \sum_{r=0}^{n-1} l_r(E).$$
\end{para}

\begin{prop} (c.f. \cite{TriCa} Proposition 3.5.1) \label{pro bun}
\\
In the notation above, the morphism
$$l(E):\mgm (\mathbb{P}(E)) \to \bigoplus_{r=0}^{n-1} \mgm(X)((r))$$
is an isomorphism.
\end{prop}

\begin{para}
Let $X$ be a smooth scheme and $E/X$ a vector bundle of rank $n$.
Put $\hat{E}:=E \times_X \mathbb{A}_X^1$. Then we have the canonical isomorphisms
$$\mgm (ThE) \overset{1}{\to} \m_X (\mathbb{P}(\hat{E})) \overset{2}{\to}
\mgm (\mathbb{P}(\hat{E})/\mathbb{P}(E)).$$
The first morphism is induced from an open immersion $E \to \mathbb{P}(\hat{E})$ which is an isomorphism by \textbf{(Exc)}. 
The second morphism is induced from the projection 
$\mathbb{P}(\hat{E})-X \to \mathbb{P}(E)$ which is an isomorphism by 
\textbf{(MV)} and \textbf{(Htp)}. From this isomorphisms, we get the following
distinguished triangle
$$\mgm(\mathbb{P}(E)) \to  \mgm (\mathbb{P}(\hat{E})) \overset{\pi_E}{\to}
\mgm (ThE) \overset{+1}{\to}.$$
Using this distinguished triangle and Proposition \ref{pro bun}, we get the
following isomorphism
$$\mgm (X)((n)) \to \bigoplus_{r=0}^n \mgm(X)((r)) \overset{l(\hat{E})^{-1}}{\to} \mgm (\mathbb{P}(\hat{E})) \overset{\pi_E}{\to} \mgm (ThE).$$
We call the inverse of this isomorphism the Thom isomorphism and denote it by $\theta(E)$.
\end{para}

\subsection{Normal cone deformation}

\begin{para}
Let $(X,Z)$ be a smooth pair of pure codimension $c$ over $k$ such that dimension of $X$ is $n$. 
We will write $B_ZX$ by the blow up of $X$ in $Z$.
Put $D_ZX:=B_{0 \times Z}(\mathbb{A}^1 \times X) - B_ZX$.
There are canonical isomorphisms 
$$\m_Z(X) \to \m_{\mathbb{A}^1_Z}(D_ZX) \leftarrow \mgm (ThN_ZX).$$ 
Hence we get the isomorphism $\m_Z(X) \to \mgm (ThN_ZX)$.
\end{para}

\subsection{Gysin triangles}

\begin{para} 
Let $(X,Z)$ be a smooth pair of pure codimension $c$ over $k$ and $i:Z \to X$ a closed immersion. 
In \cite{IntMo}, D\`eglise constructs the following functorial Gysin triangle in $\dmgm(k)$
$$\mgm(X-Z) \to \mgm(X) \overset{\scriptstyle{i^\ast}}{\to} \mgm(Z)((c)) 
\overset{\scriptstyle{\partial_{X,Z}[1]}}{\to} \mgm(X-Z)[1].$$
This is constructed from the following distinguished triangle
$$\mgm(X-Z) \to \mgm(X) \to \m_Z(X) \to \mgm(X-Z)[1]$$
and the following isomorphisms 
$$\m_Z(X) \overset{\sim}{\to} \mgm(Th(N_Z(X))) 
\overset{\scriptstyle{\theta(N_Z(X))}}{\to}  
\mgm(Z)((c))$$
where the first isomorphism is induced from the normal cone deformation. 
\end{para}

\begin{para}
Let $(f,g):(Y,T) \to (X,Y)$ be a morphism of closed pairs. We assume
$Z$ (resp. $T$) is connected and smooth over $k$ of codimension $n$
in $X$ (resp. $m$ in $Y$).\\
Then we define Gysin morphism associated to $(f,g)$, denote $(f,g)_!$ by
the following commutative diagram:
$$\xymatrix{ 
\m_T(Y) \ar[r]^{\scriptstyle{(1)}}_{\sim} \ar[d]_{\scriptstyle{(f,g)_\ast}} & \mgm(Th(N_TY))  
\ar[r]^{\scriptstyle{\theta(N_TY)}}_{\sim} \ar[d]&
\mgm(T)((m)) \ar[d]^{\scriptstyle{(f,g)_!}}\\
\m_Z(X) \ar[r]^{\scriptstyle{(1)}}_{\sim} & \mgm(Th(N_ZX))  
\ar[r]^{\scriptstyle{\theta(N_ZX)}}_{\sim} &
\mgm(Z)((n))
}$$
where morphisms (1) are isomorphisms induced from the morphisms of normal
cone deformations.
\end{para}

\begin{para} \label{Gysin covariance}
In the notation above, consider $i:Z \to X$, $j:X-Z \to X$,
$k:T \to Y$ and $l:Y-T \to T$ the canonical immersions. 
The following diagram is commutative:
$$\xymatrix{
\mgm(Y-T) \ar[r]^{\scriptstyle{\mgm(l)}} \ar[d]_{\scriptstyle{\mgm(h)}} & \mgm(Y) \ar[d]^{\scriptstyle{\mgm(f)}} 
\ar[r]^{\scriptstyle{k^\ast}} &
\mgm(T)((m)) \ar[r]^{\scriptstyle{\partial_{Y,T}[1]}} 
\ar[d]_{\scriptstyle{(f,g)_!}} & \mgm(Y-T)[1] \ar[d]^{\scriptstyle{\mgm(h)[1]}}\\
\mgm(X-Z) \ar[r]_{\scriptstyle{\mgm(j)}} & \mgm(X) 
\ar[r]_{\scriptstyle{i^\ast}} &
\mgm(Z)((n)) \ar[r]_{\scriptstyle{\partial_{X,Z}[1]}} & \mgm(X-Z)[1]
}$$
\end{para}

\begin{lem} \label{Gysin example}
Let $x:\Spec k \to \a1 \overset{j}{\hookrightarrow} \mathbb{P}^1$ be a $k$-rational point. Then the following diagram is commutative.
$$\xymatrix{
\mgm(\mathbb{P}^1) \ar[r]^{\scriptstyle{x^\ast}} \ar[d]_{\scriptstyle{\id}} & \mgm(\Spec k)((1)) \ar[d]^{\wr}\\
\mgmc(\mathbb{P}^1) \ar[r]_{\ \ \ \ \scriptstyle{j^\ast}} & \mgmc(\a1) 
}$$ 
where the vertical isomorphism $\mgm(\Spec k)((1)) \overset{\sim}{\to} \mgmc(\a1)$ is defined in \ref{Tate object}.
\end{lem}

\begin{proof}
It is just a matter of considering two split distinguished triangles below
$$\xymatrix{
\mgm(\a1) \ar[r] \ar[d]_{\scriptstyle{\mgm(p)}} \ar@{}[dr] |{1}
& \mgm(\mathbb{P}^1_{k}) \ar[r]^{\scriptstyle{x^\ast}} \ar[d]^{\scriptstyle{\id}} & \mgm(\Spec k)((1)) \ar[r] \ar[d] & \mgm(\a1)[1] \ar[d]^{\scriptstyle{\mgm(p)[1]}}\\
\mgm(\Spec k) \ar[r]_{\scriptstyle{\mgm(x)}} & \mgm(\mathbb{P}^1_{k}) \ar[r]_{\ \ \ \scriptstyle{j^\ast}} & \mgmc(\a1) \ar[r] & \mgm(\Spec k)[1]
}$$
where the commutativity of 1 follows from Example \ref{indep rat pt}.
\end{proof}

\begin{para} \label{Gysin contravariance}
In the notation above, if $f$ is finite equidimensional,
then we define $(f,g)^!$ by the following commutative diagram:
$$\xymatrix{
\m_T(Y) \ar[r]^{\scriptstyle{(1)}}_{\sim} \ar[d]_{\scriptstyle{(f,g)^\ast}} & \mgm(Th(N_TY))  
\ar[r]^{\scriptstyle{\theta(N_TY)}}_{\sim} \ar[d] &
\mgm(T)((m)) \ar[d]^{\scriptstyle{(f,g)^!}}\\
\m_Z(X) \ar[r]^{\scriptstyle{(1)}}_{\sim} & \mgm(Th(N_ZX))  
\ar[r]^{\scriptstyle{\theta(N_ZX)}}_{\sim} &
\mgm(Z)((n))
}$$
where the morphisms (1) are the isomorphisms induced from the morphisms of normal cone deformations.\\
The following diagram is commutative:
$$\xymatrix{
\mgm(Y-T) \ar[r]^{\scriptstyle{\mgm(l)}} \ar[d]_{\scriptstyle{\mgm({}^t h)}}
& \mgm(Y) \ar[r]^{\scriptstyle{k^\ast}} \ar[d]^{\scriptstyle{\mgm({}^t f)}}
& \mgm(T)((m)) \ar[r]^{\scriptstyle{\partial_{Y,T}[1]}} \ar[d]_{\scriptstyle{(f,g)^!}} 
& \mgm(Y-T)[1] \ar[d]^{\scriptstyle{\mgm(h)[1]}}\\
\mgm(X-Z) \ar[r]_{\scriptstyle{\mgm(j)}} & \mgm(X) 
\ar[r]_{\scriptstyle{i^\ast}} &
\mgm(Z)((n)) \ar[r]_{\scriptstyle{\partial_{X,Z}[1]}} & \mgm(X-Z)[1]
}$$
\end{para}

\begin{prop}(\cite{MotGe} Proposition 2.5.2) \label{Gys contra2} $ $
\\
In the notation above, if $(f,g)$ is Cartesian, then 
$$(f,g)^!=\mgm({}^t g)((m)).$$
\end{prop}

\noindent
Next we cite the some proposition in \cite{MotGe}. This is needed to prove
that the motivic reciprocity law 
implies the Weil reciprocity law for Milnor $K$-groups. 

\begin{para} \label{rho1}
By \ref{Tate object}, we have the following distinguished triangle
$$\mathbb{Z} \to \mgm(\G_m) \overset{\rho_1}{\to} \mathbb{Z}\{1\} \overset{+1}{\to}$$ 
where $\mathbb{Z} \to \mgm(\G_m)$ is induced from the unit morphism
$\Spec k \overset{1}{\to} \G_m$.
\end{para}

\begin{prop} (c.f. \cite{MotGe} Proposition 2.6.6) \label{specialization} 
\\
Let $(X,Z)$ be a smooth closed pair of codimension $1$. We denote
$i:Z \to X$ a closed immersion and $j:X-Z \to X$ a canonical open immersion.\\
Suppose there is a regular function $\pi:X \to \mathbb{A}^1_k$ which parameterizes $Z$. Hence we have a morphism $\pi=\pi|_{X-Z}:X-Z \to \G_m$.
Then the following diagram is commutative
$$\xymatrix{
\mgm(Z)\{1\} \ar[r]^{\scriptstyle{\partial_{X,Z}}} \ar[d]_{\scriptstyle{\mgm(i)\{1\}}} 
& \mgm(X-Z) \ar[d]^{\scriptstyle{\mgm(\pi) \circ \rho_1 \boxtimes \id_{\mgm(X-Z)}}}\\
\mgm(X)\{1\} & \mgm(X-Z)\{1\} \ar[l]^{\scriptstyle{\mgm(j)\{1\}}}.
}$$  
\end{prop}

\section{Motivic cohomology groups attached to pointed smooth curves}

\subsection{Definition}

\begin{para}
For pointed smooth curves $(C_1,x_1)$,\ldots,$(C_r,x_r)$ over field $k$, 
we define a motivic complex 
$\mathbb{Z}((C_1,x_1) \wedge \ldots \wedge (C_r,x_r))$, or
$\mathbb{Z}(C_1 \wedge \ldots \wedge C_r)$ in $\dmeff(k)$ as follows
$$ \mathbb{Z}(C_1 \wedge \ldots \wedge C_r)=C^{\ast}(\ztr(C_1,x_1) \otimes \ldots \otimes \ztr(C_r,x_r))[-r]$$ 
\end{para}

\begin{para}
The restriction $\mathbb{Z}(\underset{i=1}{\overset{r}{\wedge}} (C_i,x_i))|_X$ of 
$\mathbb{Z}(\underset{i=1}{\overset{r}{\wedge}} (C_i,x_i))$ to the Zariski site of $X \in Sm/k$ is a complex of sheaves in Zariski topology and the motivic cohomology groups 
$\Homo_\mathcal{M}^n(X,\underset{i=1}{\overset{r}{\wedge}}(C_i,x_i))$, or
$\Homo_\mathcal{M}^n(X,\underset{i=1}{\overset{r}{\wedge}} C_r)$ 
are defined to be the hyper cohomology of the motivic complexes 
$\mathbb{Z}(\underset{i=1}{\overset{r}{\wedge}} (C_i,x_i))$ with respect to 
Zariski topology:
$$\Homo_\mathcal{M}^n(X,\underset{i=1}{\overset{r}{\wedge}} (C_i,x_i))=
\mathbb{H}^n_{Zar}(X,\mathbb{Z}(\underset{i=1}{\overset{r}{\wedge}} (C_i,x_i))|_X).$$
\end{para}

\begin{para} \label{com}
As in \cite{TriCa} we have
$$\mathbb{H}^n_{Nis}(X,\mathbb{Z}(\underset{i=1}{\overset{r}{\wedge}} (C_i,x_i))|_X)=
\Hom_{\dmeff(k)}(\mgm(X),\mathbb{Z}(\underset{i=1}{\overset{r}{\wedge}} (C_i,x_i))[n])$$
Moreover if $k$ is a perfect field  as in \cite{CohTh} Proposition 3.1.11, 
we have
$$\mathbb{H}^n_{Zar}(X,\mathbb{Z}(\underset{i=1}{\overset{r}{\wedge}} (C_i,x_i))|_X)=
\mathbb{H}^n_{Nis}(X,\mathbb{Z}(\underset{i=1}{\overset{r}{\wedge}} (C_i,x_i))|_X).$$
\end{para}

\subsection{Fundamental properties}

\begin{para}(Product structure) \label{product}
Let $(C_1,a_1), \ldots, (C_s,a_s), (D_1,b_1),\ldots,(D_t,b_t)$ be pointed smooth curves over field $k$.
As usual motivic complexes, we have canonical morphism
$$\mathbb{Z}(\underset{i=1}{\overset{s}{\wedge}} (C_i,a_i)) \otimes \mathbb{Z}(
\underset{j=1}{\overset{t}{\wedge}}
(D_j,b_j)) \to 
\mathbb{Z}(\underset{i=1}{\overset{s}{\wedge}} (C_i,a_i) \wedge 
\underset{j=1}{\overset{t}{\wedge}} (D_j,b_j)).$$
Hence we get for any $X \in Sm/k$ the pairing
$$\Homo^p_\mathcal{M}(X, \underset{i=1}{\overset{s}{\wedge}}(C_i,a_i)) \otimes 
\Homo^q_\mathcal{M}(X,\underset{j=1}{\overset{t}{\wedge}}  (D_j,b_j))
\to \Homo^{p+q}_\mathcal{M}(X,\underset{i=1}{\overset{s}{\wedge}} (C_i,a_i) \wedge 
\underset{j=1}{\overset{t}{\wedge}} (D_j,b_j)).$$
\end{para}

\begin{para} \label{for field}
Let $(C_1,a_1),\ldots,(C_r,a_r)$ be pointed smooth curves over field $k$. For any field extension $L/k$, we abbreviate 
$\Homo^p_\mathcal{M}(\Spec L,\underset{i=1}{\overset{r}{\wedge}} 
(C_i \times_k L,a_i \times_k \id_L )$ as
$\Homo^p_\mathcal{M}(L,\underset{i=1}{\overset{r}{\wedge}} C_i)$. By definition, we have
$$\Homo^p_\mathcal{M}(L,\underset{i=1}{\overset{r}{\wedge}} C_i)=
\Homo_{i-p}(C_\ast(\bigotimes_{i=1}^r \ztr(C_i \times_k L,a_i \times_k \id_L))(\Spec L)).$$  
\end{para}

\begin{para}(Norm map) \label{norm map}
In \ref{for field}, if we assume $L/k$ is a finite field extension, then using the description of \ref{for field} and the proper push-forward of cycles induces a map 
$$\Norm_{L/k}:\Homo^p_\mathcal{M}(L,\underset{i=1}{\overset{r}{\wedge}} C_i) \to 
\Homo^p_\mathcal{M}(k,\underset{i=1}{\overset{r}{\wedge}} C_i).$$ 
From the corresponding properties of proper push-forward, the 
following properties  
are immediately verified.\\
For finite field extension $k \subset L \subset M$ and 
$x \in \Homo^p_\mathcal{M}(M,\underset{i=1}{\overset{r}{\wedge}} C_i)$ and 
$y \in \Homo^p_\mathcal{M}(L,\underset{i=1}{\overset{r}{\wedge}} C_i)$ then
we have\\
(1) $\Norm_{M/L}(y_M \centerdot x)=y \centerdot \Norm_{M/L}(x)$ and  
$\Norm_{M/L}(x \centerdot y_M)=\Norm_{M/L}(x) \centerdot y$\\
(2) $\Norm_{M/k}(x)=\Norm_{M/L}(\Norm_{L/k}(x))$\\
(3) If $M/k$ is a normal extension, we have
$$\Norm_{L/k}(x)_M=[L:k]_{insep}\underset{j:M \hookrightarrow L}{\Sigma} j^\ast(x).$$ 
\end{para}

\begin{ex} \label{weight one}
Let $(C,x)$ be a pointed projective smooth curve over $k$. Then we have 
$$\Homo^1_\mathcal{M}(k,(C,x))=\Ker(\ch_0 (C) \overset{\text{deg}}{\to} 
\mathbb{Z}).$$
\end{ex}

\begin{ex} \label{weight one 2}
Let $X$ be a smooth curve over $k$. A good compactification of $X$ is a pair $(\bar{X},X_\infty)$ such that
there is an open embedding 
$X \overset{j}{\hookrightarrow} \bar{X}$, $\bar{X}$ is proper non-singular curve over $k$ and $X_\infty=\bar{X}-X$ has an affine open neighborhood in 
$\bar{X}$.\\
Let $(C,x)$ be a pointed smooth affine curve with a good compactification
$(\bar{X},X_\infty)$, then
$$\Homo^1_\mathcal{M}(k,(C,x))=\Ker(\Pic (\bar{X},X_\infty) \overset{\text{deg}}{\to} 
\mathbb{Z})$$
where $\Pic (\bar{X},X_\infty)$ is the relative Picard group. The elements of
$\Pic (\bar{X},X_\infty)$ are the isomorphism classes $(\mathcal{L},t)$ of line bundle
$\mathcal{L}$ on $\bar{X}$ with a trivialization $t$ on $X_\infty$. 
\end{ex}

\begin{para} \label{norm ccompatibility}
In \ref{weight one}, for finite field extension $L/k$, using the property \ref{norm map} (3), we know through isomorphisms in \ref{weight one}, norm maps in \ref{norm map} and classical one are compatible. 
\end{para}

\section{Calculation of motivic cohomology groups attached to pointed smooth curves}

\subsection{Pro-motives}

In this subsection, we will briefly review the result of \cite{MotGe}.

\begin{para}
Let $\mathcal{A}$ be a tensor triangulated category. 
We consider $\proa$ the pro-category of $\mathcal{A}$.
Then the following facts are fundamental.\\
1. $\proa$ is additive.\\
2. The shift functor of $\mathcal{A}$ induces an auto-functor of $\proa$.\\
3. There is a unique tensor structure over $\proa$ such that 
$\otimes$ commutes projective limits.
\end{para}

\begin{para}
In the notation above, we call any triangle in $\proa$ isomorphic to
formal projective limit of distinguished triangles of $\mathcal{A}$  a pro-distinguished triangle. Let  
$H:\mathcal{A}^{op} \to \mathcal{A}b$ be a cohomological functor. 
Then the functor
$$\overline{H}:(\proa)^{op} \ni (X_i)_{i \in I} \mapsto
\injlim_{i \in I^{op}}H(X_i) \in \mathcal{A}b$$
sends pro-distinguished triangles to long exact sequences. 
\end{para}

\begin{para}
Let $\mathcal{O}$ be a $k$-algebra. We say $\mathcal{O}$ is local smooth over
$k$ iff there is an formally smooth of finite type $k$-algebra $A$, and a prime
ideal $x$ of $A$ and an isomorphism $\mathcal{O} \overset{\sim}{\to} A_x$.\\
Since $k$ is perfect, $\mathcal{O}$ is local smooth iff it is regular and
essentially of finite type. 
\end{para}

\begin{para}
Let $\mathcal{O}$ be a local smooth $k$-algebra. A model of $\mathcal{O}/k$
is a pair $(X,x)$ consist of  a smooth scheme $X$ and a morphism 
$x:\Spec \mathcal{O} \to X$ such that, if we write the image of closed point of
$\Spec \mathcal{O}$ to $x$, induced morphism 
$x^{\sharp}:\mathcal{O}_{X,x} \to \mathcal{O}$ is an isomorphism.\\
We put
$$\mathcal{M}^{\text{smo}}(\mathcal{O}/k):=\{
A \subset \mathcal{O};(\Spec A, 
\ \Spec \mathcal{O} \overset{\overset{\text{obvious}}{\text{map}}}{\to} \Spec A)\text{ is a model of}\ \mathcal{O}/k\}.$$   
$\mathcal{M}^{\text{smo}}(\mathcal{O}/k)$ is not empty and filtrant for inclusion.
(c.f. \cite{MotGe} Lemma 3.1.5).
\end{para}

\begin{para}
1. Let $\mathcal{O}$ be a local smooth $k$-algebra. We consider a pro-object
of $\sm/k$
$$(\mathcal{O}):=\{\Spec A\}_{A \in \mathcal{M}^{\text{smo}}(\mathcal{O}/k)}.$$
2. Let $X$ be a smooth scheme and $x \in X$, we define localization of $X$
in $x$ as a pro-object of $\sm/k$. 
$$ X_x:=\{U\}_{x \in U \subset X}$$ 
where $U$ runs through the open neighborhood of $x$.\\
\end{para}

\begin{para} (c.f. \cite{MotGe} Lemma 3.1.8)
Let $\mathcal{O}$ be a local smooth $k$-algebra, and $(X,x)$ a model of 
$\mathcal{O}$. Then $x$ induces a canonical isomorphism
$$(\mathcal{O}) \to X_x.$$
\end{para}

\begin{para}
Let $\mathcal{O}$ be a local smooth $k$-algebra and $n,m \in \mathbb{Z}$.
We consider a pro-object of $\prodmgm(k)$
$$\mgm(\Spec \mathcal{O})(n)[m]:=
\{\mgm(\Spec A)(n)[m]\}_{A \in \mathcal{M}^{\text{smo}}
(\mathcal{O}/k)}.$$ 
\end{para}

\noindent
Next we define a residue morphism associated to a discrete valuation.

\begin{para}
Let $E/k$ be a field extension of finite type, $v$ a valuation of $E/k$,
$\mathcal{O}_v$ a valuation ring of $v$ and $(X,t)$ a $k$-model of 
$\mathcal{O}_v$.
We say a special point of $(X,t)$ for image of closed point of
$\Spec \mathcal{O}_v$ for $t$ and denote by $s$.\\
We say that $(X,t)$ is a strict $k$-model of $\mathcal{O}_v$ iff closure 
$\bar{\{s\}}$ in $X$ is a smooth scheme.\\
Any discrete valuation ring $\mathcal{O}_v$ essentially of finite type
over $k$ admits a strict $k$-model. (c.f. \cite{MotGe} Lemma 4.5.3).
\end{para}

\begin{para}
Let $E/k$ be a field extension of finite type, $v$ a valuation of $E/k$,
$\mathcal{O}_v$ a valuation ring of $v$ and $(X,t)$ a strict $k$-model of 
$\mathcal{O}_v$. Put $Z:=\bar{\{s\}}$. Since $(X,Z)$ is a smooth closed pair
of codimension $1$. We have a distinguished triangle of the form
$$\mgm(Z)\{1\} \overset{\partial_{X,Z}}{\to} \mgm(X-Z) \overset{j_\ast}{\to}
\mgm(X) \overset{+1}{\to}.$$
Since this triangle is natural for inclusions of open sets in $X$. 
(c.f. \ref{Gysin covariance})
Considering
a cofiltrant system of open neighborhoods of $s$ in $X$, we get a 
pro-distinguished triangle
$$\mgm(Z_s)\{1\} \overset{\partial_{X_s,Z_s}}{\to} \mgm(X_s-Z_s)
\overset{j_{s \ast}}{\to} \mgm(X_s) \overset{+1}{\to}.$$
Since $(X,s)$ is a $k$-model of $\mathcal{O}_v$, a morphism 
$s:\Spec \mathcal{O}_v \to X$ induces an isomorphism of pro-object
$(\mathcal{O}_v) \to X_s$; Hence we get a pro-distinguished triangle
isomorphic to the form
$$\mgm(\Spec k(v))\{1\} \overset{\partial_{X_s,Z_s}}{\to} \mgm(\Spec E)
\overset{i^\sharp}{\to} \mgm(\Spec \mathcal{O}_v) \overset{+1}{\to}$$
where $E$ (resp. $k(v)$) is a fraction field (resp. residue field) of $v$,
$i:\mathcal{O}_v \to E$ is a canonical inclusion.
\end{para}

\begin{lem} (c.f. \cite{MotGe} Lemma 4.5.5.) \label{independence} \\
Let $E/k$ be a field extension of finite type, $v$ a discrete
valuation of $E/k$,
$\mathcal{O}_v$ a valuation ring of $v$ and $k(v)$ is this residue field.\\
Then adopting the notation above, if $(X,s)$ and $(Y,t)$ are two strict
$k$-models of $\mathcal{O}_v$, put $Z:=\bar{\{s\}}$ closure in $X$ and 
$T:=\bar{\{t\}}$ closure in $Y$. Then we have
$$\partial_{X_s,Z_s}=\partial_{Y_t,T_t}.$$ 
\end{lem}

\begin{para} \label{residuedef}
Let $E/k$ be a field extension of finite type, $v$ a discrete 
valuation of $E/k$.\\
We define a residue morphism associated to $v$, denoted by $\partial_v$
defined by $\partial_v:=\partial_{(X_s,Z_s)}$ where $(X,s)$ is a strict $k$-model of valuation ring of $v$. 
By Lemma \ref{independence} this does not depend on a choice of a strict
$k$-model of valuation ring of $v$.\\
So we have the following pro-distinguished triangle of the form
$$\mgm(\Spec k(v))\{1\} \overset{\partial_v}{\to} \mgm(\Spec E)
\overset{i^\sharp}{\to} \mgm(\Spec \mathcal{O}_v) \overset{+1}{\to}.$$
\end{para}

\noindent
Having defined residue morphisms, we explain the connection of Milnor $K$-groups and Hom sets in $\prodmgm(k)$.

\begin{para}
Using the distinguished triangle in \ref{rho1} and the definition of 
tensor structure in $\dmeff(k)$,
for any $n \in \mathbb{N}$, there is a distinguished triangle of the form
$$\bigoplus_{i=1}^n \mgm(\G_m^{n-1}) \to \mgm(\G_m^{n})
\overset{\rho_n}{\to} \mathbb{Z}\{n\} \overset{\scriptstyle{+1}}{\to}.$$ 
where the first morphism is induced from sum of $n$ closed immersions
$$\iota_i:=\id \times \id \times \ldots \times \overset{i}{1} \times \ldots
\times \id \times \id.$$  
\end{para}

\begin{para} \label{external product}
Let $E/k$ be a field extension of finite type. 
$f:\mgm(\Spec E) \to \mathcal{M}$ and $g:\mgm(\Spec E) \to \mathcal{N}$
are morphisms in $\prodmgm(k)$. Then we can extend the definition of
external cup product $f \boxtimes g: \mgm(\Spec E) \to \mathcal{M}
\otimes \mathcal{N}$.\\
If $\mathcal{M}=\mathbb{Z}\{p\}$ and $\mathcal{N}=\mathbb{Z}\{q\}$, we have
a canonical isomorphism $\mathbb{Z}\{p\} \otimes \mathbb{Z}\{q\} 
\overset{\sim}{\to} \mathbb{Z}\{p+q\}$. Then we have the following identity
$$\alpha \boxtimes \beta = - \beta \boxtimes \alpha.$$
(c.f. \cite{MotGe} Remarque 4.4.2.)  
\end{para}

\begin{para}
Let $E/k$ be a field extension of finite type. Then we have a morphism
\begin{eqnarray*}
(E^\times)^n  & \overset{\sim}{\to} &
\Hom(\Spec E, \G_m^n)\\
& \overset{\scriptstyle{\mgm}}{\to} &
\Hom_{\prodmgm(k)}(\mgm(\Spec E),\mgm(\G_m^n))\\ 
& \overset{\scriptstyle{\Hom(\mgm(\Spec E),\rho_n)}}{\longrightarrow} &
\Hom_{\prodmgm(k)}(\mgm(\Spec E),\mathbb{Z}\{n\}).
\end{eqnarray*}
This map induces a morphism
$$\alpha:K_n^M(E) \to \Hom_{\prodmgm(k)}(\mgm(\Spec E),\mathbb{Z}\{n\}).$$
\end{para}

\begin{para}
In the notation above, for any $x \in (E^\times)^{\otimes n}$ and
$y \in (E^\times)^{\otimes m}$, we have
$$\alpha(x \otimes y)= \alpha(x) \boxtimes \alpha(y).$$ 
\end{para}

\begin{thm} (c.f. \cite{MotGe} Theorem 4.4.4)  \label{Milkisom} $ $
\\
(Need the assumption of perfectness of $k$.) In the notation above,
$$\alpha:K_{\ast}^M(E) \to 
\Hom_{\prodmgm(k)}(\mgm(\Spec E),\mathbb{Z}\{\ast\})$$ 
is an algebra isomorphism.
\end{thm}

\subsection{Motivic reciprocity law}

The classical theorems ``Weil reciprocity law'' and ``residue formula'' are unified using Milnor K-groups. More precisely, the following statement is known. (c.f. \cite{Sus82})

\begin{parathm}(Reciprocity law for Milnor $K$-groups)
\\
Let $K$ be an algebraic function field over a field $k$. Then the following 
composition are the zero maps for all non-negative integers $n$. 
$$K^M_{n+1}(K) \overset{\bigoplus \partial_v}{\to} 
\underset{v}{\bigoplus} 
K^M_n(k(v)) \overset{\Sigma \Norm_{k(v)/k}}{\to} K^M_n(k)$$
\end{parathm} 

\noindent 
In this subsection, we will prove more fundamental style of the following reciprocity law.

\begin{thm}(Motivic reciprocity law) \label{mot rec law}
\\
The following composition
$$\mgm(\Spec k)\{1\} \overset{\scriptstyle{\Sigma \Norm_{k(v)/k}\{1\}}}
{\longrightarrow} \underset{v}{\overset{\sim}{\prod}}
\mgm(\Spec k(v))\{1\} \overset{\scriptstyle{\overset{\sim}{\prod} \partial_v}}{\longrightarrow} \mgm(\Spec K)$$
is the zero map in $\prodmgm(k)$. 
\end{thm}

\begin{para}
Let $K/k$ be a field extension of transcendental degree one. 
Let $C/k$ be a projective nonsingular curve such that $K(C)=K$. As in the previous subsection, we can construct the following pro-distinguished triangle
in $\prodmgm (k)$.
\begin{multline*}
\mgm(\Spec K) \to \mgm(C)\\
\to \underset{\scriptstyle{x \in C:\text{closed points}}}{\overset{\sim}{\prod}}\mgm(\Spec k(x))((1))
\overset{\scriptstyle{\overset{\sim}{\prod} \partial_x[1]}}{\to} \mgm(\Spec K)[1]
\end{multline*}
This is constructed as follows:
For any closed set $Z \subset C$, there is the Gysin triangle 
$$\mgm(C-Z) \to \mgm(C) \to \underset{\scriptstyle{x \in Z}}
{\bigoplus}\mgm(\Spec k(x))((1))
\overset{\scriptstyle{\partial_{C,Z}[1]}}{\to} \mgm(C-Z)[1]$$
and we consider 
$\mgm(\Spec K)=\{\mgm(C-Z)\}_{\scriptstyle{Z \subset C:\text{closed subsets}}} \in 
\prodmgm (k)$
\end{para}

\begin{lem} $ $
\\
In the notation above, for any closed point $x \in C$, the diagram of structure morphisms
$$\xymatrix{ 
{\Spec k(x)} \ar[r]^{\scriptstyle{i}} \ar[dr] & C \ar[d]_{\scriptstyle{p}}\\
& {\Spec k}\\
}$$
induces the following commutative diagram:
$$\xymatrix{ 
\mgm (\Spec k(x))((1)) & \mgm (C) \ar[l]_{\scriptstyle{i^\ast}}\\
& \mgm (\Spec k)((1)) \ar[lu]^{\scriptstyle{\Norm_{k(x)/k}((1))}} \ar[u]_{\scriptstyle{p^\ast}}.\\
}$$
\end{lem}

\begin{proof}
First choose a finite equidimensional morphism $C \overset{\pi}{\to} \mathbb{P}^1$ which is unramified at every points over $\pi(x)$. (This can be done by using Bertini theorem.) Next using \ref{fund of flat pullback}
and Proposition \ref{Gys contra2}, we may assume $C=\mathbb{P}^1$. Replacing $\mathbb{P}^1$ by $\mathbb{P}^1_{k(x)}$ and using Proposition \ref{Gys contra2} again, we may assume $k(x)=k$.
In this case, $i^\ast \circ p^\ast =\id$ by Lemma \ref{flat pillback example} and \ref{Gysin example}.
\end{proof}

\begin{para}
Hence we get the following diagram:
$$\xymatrix{
&  \mgm (\Spec k)((1)) \ar[ld]_{\scriptstyle{p^\ast}} \ar[d]^{\scriptstyle{\Sigma \Norm_{k(x)/k}((1))}}\\
\mgm (C) \ar[r] & 
\underset{\scriptstyle{x \in Z}}{\bigoplus} \mgm (\Spec k(x))((1)) \ar[r] & \mgm (C-Z)[1].\\
}$$
\noindent
Taking a limit with respect to $Z$, we get the motivic reciprocity law.
\end{para}

\noindent
Next we prove that the motivic reciprocity law implies the Weil reciprocity law
for Milnor $K$-groups.

\begin{lem} \label{some comm} $ $
\\
In the notation above, let $v$  be a valuation of $K/k$,
$\mathcal{O}_v$ a valuation ring of $v$, $\pi$ a uniformizer element
of $\mathcal{O}_v$. Then the following diagram is commutative.
$$\xymatrix{
\mgm (\Spec k(v))\{1\} \ar[r]^{\scriptstyle{{\partial_v}}} \ar[d] & \mgm (\Spec K) \ar[d]^{\scriptstyle{{\alpha(\pi) \boxtimes \id_{\mgm (\Spec K)}}}}\\
\mgm (\Spec \mathcal{O}_v)\{1\} & \mgm (\Spec K)\{1\} \ar[l].
}$$
\end{lem}

\begin{proof}
Take a strict $k$-model of $\Spec \mathcal{O}_v$. Denote it by $(X,s)$.
If we take $X$ sufficiently small, $\pi \in \mathcal{O}_v \overset{\sim}{\to}
\mathcal{O}_{X,s}$ determines a regular function $X \to \a1$ which parameterizes
$Z$. Using Proposition \ref{specialization} and considering a cofiltrant system of open neighborhoods of $s$ in $X$, we get the following commutative diagram of pro-motives
$$\xymatrix{ 
\mgm(Z_s)\{1\} \ar[r]^{\scriptstyle{\partial_{X_s,Z_s}}} \ar[d] 
& \mgm(X_s-Z_s) \ar[d]^{\scriptstyle{\mgm (\pi) \circ \rho_1 \boxtimes \id_{\mgm (X_s-Z_s)}}}\\
\mgm (X_s)\{1\} & \mgm (X_s-Z_s)\{1\}. \ar[l] 
}$$
Hence we get the result.
\end{proof}

\begin{ex} \label{coincide}
In the notation above, for any discrete valuation $v$ of $K/k$, 
there is a commutative diagram\\
$$\xymatrix{ 
\Hom (\mgm (\Spec K),\mathbb{Z}\{n+1\}) 
\ar[d]_{\scriptstyle{{\Hom (\partial_v,\mathbb{Z}\{n+1\})}}}
& K^M_{n+1}(K) \ar[l]^{\sim}
\ar[d]^{\scriptstyle{{(-1)^n \partial_v}}}\\
\Hom (\mgm (\Spec k(v)\{1\},\mathbb{Z}\{n+1\}) & K^M_n(k(v)). \ar[l]^{\sim}
}$$ 
\end{ex}

\begin{proof}
For $u_1,\ldots,u_{n+1} \in \mathcal{O}_v^\times$ and a uniformizer $\pi$, it is enough to check the following two conditions.\\
1 $\Hom (\partial_v, \mathbb{Z}\{n+1\})(\alpha(\{u_1,\ldots,u_{n+1}\}))=0$\\
2 $\Hom (\partial_v, \mathbb{Z}\{n+1\})(\alpha(\{u_1,\ldots,u_n,\pi\}))=
(-1)^n \alpha (\{\overline{u_1},\ldots,\overline{u_n}\})$\\
To prove 1: Notice that there is a pro-distinguished triangle as follows
(c.f. \ref{residuedef})\\
$$\mgm (\Spec k(v))\{1\} \overset{\scriptstyle{\partial_v}}{\to} 
\mgm (\Spec K) \to
\mgm (\Spec \mathcal{O}_v) \overset{\scriptstyle{+1}}{\to}.$$
To prove 2: Anti-commutativity of $\boxtimes$ (c.f. \ref{external product})
and Lemma \ref{some comm}.
\end{proof}

\begin{ex} \label{Gysin compatibility}
Let $K/k$ be a field extension of transcendental degree 1, $(C,x)$ be a pointed smooth curve and $v$ a place of $K/k$. There is a tame symbol
$\partial_v:\Jac C(K_v) \otimes K_v^{\times} \to \Jac C(k(v))$. Then the following diagram is commutative.
$$\xymatrix{
\Jac C(K_v) \otimes K_v^{\times} \ar[r] \ar[dd]_{\scriptstyle{\partial_v}}
& \overset{\Hom_{\dmeff(k)}(\mgm(\Spec K_v), \mathbb{Z}(C,x)[1])}{\underset{\Hom_{\dmeff(k)}(\m (\Spec K_v), \mathbb{Z}(\G_m,1)[1])} {\otimes}} \ar[d]^{\scriptstyle{\boxtimes}}\\
& \Hom_{\dmeff(k)}(\m (\Spec K_v)[1], \mathbb{Z}(C,x)(1)[3])
\ar[d]^{\scriptstyle{-\Hom(\partial_v,\mathbb{Z}(C,x)(1)[3])}} \\
\Jac C(k(v)) \ar[r] 
& \Hom_{\dmeff(k)}(\m (\Spec k(v))(1)[2], \mathbb{Z}(C,x)(1)[3]) 
}$$   
This is proved in the same way as Example \ref{coincide}. 
\end{ex}

\begin{cor} \label{mot->Weil} $ $
\\
The motivic reciprocity law implies the Weil reciprocity law for Milnor $K$-groups.
\end{cor}

\begin{proof}
Take $\Hom(?,\mathbb{Z}\{n+1\})$ and use Theorem \ref{can thm}, 
Theorem \ref{Milkisom} and
notice Example \ref{norm} and Example \ref{coincide}. 
\end{proof}

\subsection{Main result}

\begin{para}
In this section, let $k$ be a perfect field which admits resolution of 
singularities and $(C_1,a_1),\ldots,(C_n,a_n)$ pointed projective smooth curves over $k$. 
\end{para}

\begin{para} \label{weil reciprocity 1}
Let $p:Z \to \a1$ be a finite surjective morphism and suppose that $Z$ is integral. Let $f_i \in \Hom(Z,C_i)$ and
$$p^{-1}({j})=\coprod n_i^jz_i^j \ \ \ (j=0,1)$$
where $n_i^j$ are the multiplicities of points $z_i^j=\Spec L_i^j$. Define:
$$ \phi_j=\Sigma n_i^j\{f_1,\ldots,f_n\}_{L^j_i/k}$$
then we have $$\phi_0=\phi_1$$ in $K(k,\Jac C_1,\ldots,\Jac C_n)$.\\
The proof is similar to \cite{MVW02}, p.45 Corollary 5.5. 
\end{para}

\begin{para}
As $\bigotimes_{i=1}^n \ztr (C_i,a_i)(\Spec k))$ is a quotient of the free abelian groups generated by the closed points of $C_1 \times \ldots \times C_n$ modulo the subgroup generated by all points of the form $(x_1,\ldots,a_i,\ldots,x_n)$ where the $a_i$'s can be any position. If $x$ is a closed point of $C_1 \times \ldots \times C_n$ with residue field $L$ then $x$ is defined by a canonical sequence $(x_1,\ldots,x_n) \in \Jac C_1(L) \times \ldots \times \Jac C_n(L)$.\\
Since 
$$\Homo^n_\mathcal{M}(k,\underset{i=1}{\overset{n}{\wedge}} C_i)=
\coker(\bigotimes_{i=1}^n \ztr (C_i,a_i)(\a1)
\overset{\partial_0-\partial_1}{\to}
\bigotimes_{i=1}^n \ztr (C_i,a_i)(\Spec k))$$ 
using \ref{weil reciprocity 1}, we have a natural map 
$\Homo^n_\mathcal{M}(k,\underset{i=1}{\overset{n}{\wedge}} C_i) 
\to K(k,\Jac C_1,\ldots,\Jac C_n)$.
\end{para}

\begin{para}
Using \ref{weight one}, for every finite field extension $L/k$ we have an isomorphism
$$\bigotimes_{i=1}^n \Jac C_i(L) \overset{\sim}{\to} 
\bigotimes_{i=1}^n \Homo^1_\mathcal{M}(L,C_i).$$
Combining a natural pairing in \ref{product}
$$\bigotimes_{i=1}^n \Homo^1_\mathcal{M}(L,C_i) 
\to \Homo ^n_\mathcal{M}(L,\underset{i=1}{\overset{n}{\wedge}} C_i)$$
and a norm map (c.f \ref{norm map})
$$\Norm_{L/k}:\Homo ^n_\mathcal{M}(L,\underset{i=1}{\overset{n}{\wedge}} C_i) 
\to \Homo ^n_\mathcal{M}(k,\underset{i=1}{\overset{n}{\wedge}} C_i)$$
we get a canonical map $$\underset{L/k:\text{finite extension}}{\oplus}
\bigotimes_{i=1}^n\Jac C_i(L) \to 
\Homo ^n_\mathcal{M}(k,\underset{i=1}{\overset{n}{\wedge}} C_i).$$
If we use \ref{norm map} (1), Theorem \ref{norm ccompatibility}, Theorem \ref{mot rec law} and Example \ref{Gysin compatibility}, this map should factor through the map 
$K(k,\Jac C_1,\ldots,\Jac C_n) \to \Homo ^n_\mathcal{M}(k,\underset{i=1}{\overset{n}{\wedge}} C_i).$
\end{para}

\begin{para}
Obviously the morphisms above are inverse to each other. 
Hence we get the following result.
\end{para}

\begin{thm} (Somekawa conjecture for Jacobian varieties) $ $
\\
Let $(C_1,a_1),\ldots,(C_n,a_n)$ be pointed projective smooth curves over perfect field $k$ which
admits resolution of singularities.  Then
$$K(k,\Jac C_1,\ldots,\Jac C_n) \overset{\sim}{\to} \Hom_{\dmeff(k)}
(\m(\Spec k), \mathbb{Z}(\underset{i=1}{\overset{n}{\wedge}} C_i)[n]).$$
\end{thm}

\noindent
\textbf{Acknowledgment} 
The author is greatful for Professors Takeshi Saito, Shuji Saito, Kazuya Kato, Takao Yamazaki and Kenichiro Kimura.

\begin{center}{\it
Satoshi Mochizuki\\
Graduate school of Mathematical Sciences,
The University of Tokyo, 3-8-1 Komaba, Meguro-ku
Tokyo 153-8914, JAPAN}
\\
\tt{E-mail:mochi@ms.u-tokyo.ac.jp}
\end{center}

\end{document}